\documentclass[12pt,a4paper]{article}
\usepackage{amssymb,amsfonts,amsthm,amsmath}
\usepackage{graphicx}
\usepackage{placeins}

\def\hang{\hangindent\parindent}
\def\rf{\par\noindent\hang}
\def\ssk{\smallskip}

\def\rf{\par\noindent\hang}
\def\nin{\noindent}

\newtheorem{theorem}{Theorem}

\newtheorem{corollary}{Corollary}

\begin{document}

\baselineskip=20pt

\begin{center} {\Large{\textbf {The coverage probabililty of confidence intervals in regression
after  a preliminary F test}}}
\end{center}

\bigskip
\bigskip

\noindent {
{\normalsize {\bf PAUL KABAILA$^{1*}$}}  {\footnotesize AND} {\normalsize {\bf DAVIDE FARCHIONE$^2$}}}

\bigskip

\noindent {\it 1. Department of Mathematics and Statistics, La Trobe University, Victoria 3086, Australia}

\bigskip

\noindent {\it 2.  School of Mathematics and Statistics,
University of South Australia, South Australia 5095, Australia}

\vspace{13cm}

\noindent * Author to whom correspondence should be addressed.

\noindent E-mail: P.Kabaila@latrobe.edu.au, Facsimile: 61 3 9479 2466,
Telephone: 61 3 9479 2594

\newpage



\begin{center}
 {\bf Abstract}
\end{center}

\ssk

\noindent Consider a linear regression model with
regression parameter $\beta=(\beta_1, \ldots, \beta_p)$ and independent normal
errors.
Suppose the parameter of interest is $\theta = a^T \beta$, where $a$ is specified.
Define the $s$-dimensional parameter vector $\tau = C^T \beta - t$, where $C$
and $t$ are specified.
Suppose that we carry out a preliminary F test of
the null hypothesis $H_0: \tau = 0$ against the alternative hypothesis $H_1: \tau \ne 0$.
It is common statistical practice to then construct
 a confidence interval for $\theta$ with
 nominal coverage $1-\alpha$, using the same data,
based on the assumption that the selected model had been given to us {\it a priori}
(as the true model).
 We call this the naive $1-\alpha$ confidence interval for $\theta$. This assumption is false and
it may lead to this confidence interval having
minimum coverage probability far
below $1-\alpha$, making it completely inadequate. Our aim is to compute this minimum
coverage probability.
It is straightforward to find
an expression for the coverage probability of this confidence interval
that is a multiple integral of dimension $s+1$.
However, we derive a new elegant and computationally-convenient formula for this coverage probability.
For $s=2$ this formula is a sum of a
 triple and a double integral and for all $s > 2$ this formula is a sum of a quadruple and a double integral.
  This makes it easy to compute the minimum coverage probability of the naive confidence interval,
 irrespective of how large $s$ is.
A very important practical application of this formula is to the {\sl analysis of covariance}.
In this context,  $\tau$ can be defined so that $H_0$
expresses the hypothesis of ``parallelism''.
Applied statisticians commonly recommend carrying out
 a preliminary F test of this hypothesis.
 We illustrate the application of our formula with a real-life
analysis of covariance data set and a preliminary F test for ``parallelism''.
We show that the naive 0.95 confidence interval has minimum coverage probability 0.0846,
showing that it is completely inadequate.

\bigskip

\bigskip

\nin{\sl Keywords:} analysis of covariance, naive confidence interval, preliminary F test, test for parallelism.

\vfil
\eject

\noindent {\large{\bf 1. Introduction}}

\medskip

\noindent Consider the linear regression model
 $Y = X \beta + \varepsilon$,
 where $Y$ is a random $n$-vector of responses, $X$ is a known $n \times p$ matrix with linearly
independent columns, $\beta$ is an unknown parameter $p$-vector and
$\varepsilon \sim N(0, \sigma^2 I_n)$ where $\sigma^2$ is an unknown positive parameter.
Suppose that the parameter of interest is $\theta = a^T \beta$ where $a$ is a given
$p$-vector ($a \ne 0$). We seek a $1-\alpha$ confidence interval for $\theta$.

Let the $s$-dimensional parameter vector $\tau$ be defined to be $C^T \beta - t$ where $C$ is a specified $p \times s$ matrix ($s < p$)
with linearly independent columns and $t$ is a specified $s$-vector.
Suppose that $a$ does not belong to the linear subspace spanned by the columns
of $C$.
Also suppose that we carry out a preliminary F test of
the null hypothesis $H_0: \tau = 0$ against the alternative hypothesis $H_1: \tau \ne 0$.
It is then common statistical practice to construct a confidence interval for $\theta$ with
 nominal coverage $1-\alpha$, using the same data, based on the assumption
 that the selected model had been given to us {\it a priori} (as the true model).
 We call this the naive $1-\alpha$ confidence interval for $\theta$.
In Section 2, we provide a convenient description of this confidence interval.
 This assumption is false and it can lead to the naive $1-\alpha$ confidence interval
 having minimum coverage probability far below $1-\alpha$, making it completely inadequate.
Our aim is to compute this minimum
coverage probability.
 For $s=1$, the preliminary F test is equivalent to a t test. The case of a single preliminary t test
has been dealt with by Kabaila and Giri (2009, Theorem 3). So, in the present paper,
we restrict attention to the case that $s > 1$.

Straightforward application of the methodology of Farchione (PhD thesis,
2009, Section 5.7) leads to
an expression for the coverage probability of the naive $1-\alpha$ confidence interval, for a
given value of an $s$-dimensional parameter vector, that is a multiple integral of dimension $s+1$.
Finding the minimum coverage probability using this formula becomes increasingly cumbersome
as $s$ increases due to both the need to (a) evaluate multiple integrals of dimension $s+1$ and
(b) the need to search for the minimum over a space of dimension $s$.

In Section 3, by a careful consideration of the geometry of the
situation, we derive a new elegant and computationally-convenient formula for the coverage probability
of this confidence interval for given parameter values.
For $s=2$ this formula is a sum of a
 triple and a double integral and for all $s > 2$ this formula is a sum of a quadruple and a double integral.
 This formula also shows that the coverage probability is a function of a two-dimensional parameter
 vector, irrespective of how large $s$ is. This makes it easy to compute the minimum coverage probability of the naive confidence interval,
 irrespective of how large $s$ is. Another important aspect of this formula is that it
can be used to delineate general categories of $a$, $C$ and $X$ for which the naive confidence interval
has poor coverage properties.

A very important practical application of this formula is to the {\sl analysis of covariance}.
In this context,  $\tau$ can be defined so that
$H_0$ expresses the null hypothesis of ``parallelism''.
In the applied statistics literature on the analysis of covariance it is commonly recommended
that a preliminary F test of the null hypothesis of ``parallelism'' be carried out.
See, for example, Kuehl (2000, p.563), Milliken and Johnson (2002, pp. 14 -- 17) and
Freund et al (2006, pp. 363 -- 368).
For an analysis of covariance, we can choose $a$ so that the parameter
$\theta$ is the difference in expected responses for two specified treatments, for the same specified values of the covariates.

In Section 4, we illustrate the application of the results of the paper with a real-life
{\sl analysis of covariance} data set and a preliminary F test for ``parallelism''.
We define
 $\theta$ to be
 $(\text{expected response to treatment 1}) - (\text{expected response to treatment 2})$,
 evaluated at the same specified value of the covariate.
We show that
the naive 0.95 confidence interval for $\theta$ has minimum coverage probability 0.0846,
for this specified value of the covariate. This shows that this
confidence interval is completely inadequate, for this specified value of the covariate.


\bigskip

\noindent {\large{\bf 2. Description of the naive confidence interval}}

\medskip

\noindent In this section we provide a convenient description of the naive $1-\alpha$ confidence
interval constructed after
the preliminary F test.
Let $\hat \beta$ denote the least squares estimator of $\beta$.
Define $R(\beta) = (Y - X \beta)^T (Y - X \beta)$. Let $m=n-p$.
Define $\hat \Sigma^2 = R(\hat \beta)/m = (Y - X \hat \beta)^T (Y - X \hat \beta)/m$.
Also, define $\hat \Theta = a^T \hat \beta$ and $\hat \tau = C^T \hat \beta - t$.
We suppose that the columns of the matrix $C$ are linearly independent.
We also suppose that $a$ does not belong to the linear subspace spanned by the
columns of $C$.
Now define the $(s+1) \times (s+1)$ matrix
\begin{equation*}
V = \left[\begin{matrix} v_{11} \quad v_{21}^T \\ v_{21} \quad V_{22} \end{matrix} \right ]
= \frac{1}{\sigma^2} E \left (\left[\begin{matrix} \hat \Theta - \theta \\ \hat \tau - \tau \end{matrix} \right ]
\left[\begin{matrix} \hat \Theta - \theta \quad (\hat \tau - \tau)^T \end{matrix} \right ]
\right ).
\end{equation*}
Note that $v_{11} = a^T (X^T X)^{-1} a$, $v_{21} = C^T (X^T X)^{-1} a$ and $V_{22} = C^T (X^T X)^{-1} C$.

Define $\beta^*$ to be the value of $\beta$ minimizing $R(\beta)$ subject to the restriction that
$\tau = C^T \beta - t = 0$. As is well known (see e.g. Graybill, 1976, p.222)
\begin{align*}
\beta^* &= \hat \beta - (X^T X)^{-1} C (C^T (X^T X)^{-1} C)^{-1} \hat \tau \\
R(\beta^*) &= R(\hat \beta) + (\hat \beta - \beta^*)^T (X^T X) (\hat \beta - \beta^*).
\end{align*}
The standard test statistic for testing $H_0 : \tau =0$ against $H_1 : \tau \ne 0$ is
\begin{equation*}
F = \frac{(\hat \beta - \beta^*)^T X^T X (\hat \beta - \beta^*) / s}{\hat \Sigma^2}
= \frac{\hat \tau^T V_{22}^{-1} \hat \tau / s}{\hat \Sigma^2}.
\end{equation*}
This test statistic has an $F_{s, m}$ distribution under $H_0$.
Suppose that we reject $H_0$ when $F > \ell$ and accept $H_0$ otherwise, where
$\ell$ is a specified positive number.

Define $\Theta^* = a^T \beta^*$.
Also define the quantile $t(m)$ by the requirement that
$P \big(-t(m) \le T \le t(m) \big) = 1-\alpha$ for $T \sim t_m$.
The naive $1-\alpha$ confidence interval for $\theta$ is obtained as follows.

Suppose that $F > \ell$. The confidence interval is constructed on the assumption that
$\tau = 0$ is not necessarily true. In this case, the naive $1-\alpha$ confidence interval
is the usual $1-\alpha$ confidence interval for $\theta$ based on fitting the full model,
\begin{equation*}
I = \big [\hat \Theta - t(m) \sqrt{v_{11}} \, \hat \Sigma, \,
\hat \Theta + t(m) \sqrt{v_{11}} \, \hat \Sigma \big].
\end{equation*}

Now suppose that $F \le \ell$. The confidence interval is constructed on the assumption that
$\tau = 0$. If $\tau = 0$ then $R(\beta^*)/\sigma^2 \sim \chi^2_{m+s}$ and
Var$(\Theta^*) = \sigma^2 \big(v_{11} - v_{21}^T V_{22}^{-1} v_{21} \big)$.
Note that $\Theta^*$ and $R(\beta^*)$ are independent random variables.
We use the notation $[a \pm b]$ for the interval $[a-b, a+b]$ ($b > 0$).
In this case, the naive $1-\alpha$ confidence interval for $\theta$ is
\begin{align}
\label{J}
J &= \left [ \Theta^* \pm t(m+s) \sqrt{\frac{R(\beta^*)}{m+s}}
\sqrt{v_{11} - v_{21}^T V_{22}^{-1} v_{21}} \right ] \notag \\
&= \left [ \Theta^* \pm t(m+s) \sqrt{\frac{R(\hat \beta)+\hat \tau^T V_{22}^{-1} \hat \tau }{m+s}}
\sqrt{v_{11} - v_{21}^T V_{22}^{-1} v_{21}} \right ].
\end{align}

\bigskip

\noindent {\large{\bf 3. The coverage probability of the naive confidence interval}}

\medskip

\noindent Define $b = v_{11}^{-1/2} V_{22}^{-1/2} v_{21}$ and
 $W = \hat \Sigma/\sigma$.
Let $f_W$ denote the probability density function of $W$.
Define $||b|| = \sqrt{b_1^2 + \cdots + b_s^2}$.
 Thus
 \begin{align*}
 ||b||^2 &= v_{11}^{-1} v_{21}^T V_{22}^{-1} v_{21} \\
 &= \frac{a^T (X^T X)^{-1} C \big (C^T (X^T X)^{-1} C \big)^{-1} C^T (X^T X)^{-1} a}{a^T (X^T X)^{-1} a}.
 \end{align*}
Since Var$(\Theta^*) = \sigma^2 \big(v_{11} - v_{21}^T V_{22}^{-1} v_{21} \big) \ge 0$,
$||b|| \in [0,1]$.
The assumption that the vector $a$ does not belong to the linear subspace spanned by the columns
of $C$ implies that $||b|| > 0$. So, we may assume that $||b|| \in (0,1]$.
Now define
\begin{equation*}
i \big(x, w; ||b|| \big) = P \big( - t(m) w + x \le Z \le t(m) w + x \big),
\end{equation*}
where $Z \sim N(0, 1 - ||b||^2)$, and
\begin{align*}
j \big(x, y, w; ||b|| \big) = P \Bigg( &x - t(m+s) \sqrt{\frac{m w^2 + y}{m+s}} \sqrt{1 - ||b||^2}
 \le Z \le \\
 &\phantom{1234567890123456789012} x + t(m+s) \sqrt{\frac{m w^2 + y}{m+s}} \sqrt{1 - ||b||^2} \Bigg).
\end{align*}
Define $f_R$ to be the probability density function of $\sqrt{R^2}$ when $R^2 \sim \chi^2_s$.
Let $B(a,b)$ denote the beta function.
Define the probability density function $f_{T_1}$ to be
\begin{equation*}
f_{T_1}(t_1) =
\begin{cases}
{\displaystyle\frac{\pi}{B(\frac{1}{2}, \frac{s-1}{2})}} \sin^{s-2}(\pi t_1) &\text{for } 0 \le t_1 \le 1 \\
0 &\text{otherwise.}
\end{cases}
\end{equation*}
For $s \ge 3$, define the probability density function $f_{T_2}$ to be
\begin{equation*}
f_{T_2}(t_2) =
\begin{cases}
{\displaystyle\frac{\pi}{B(\frac{1}{2}, \frac{s-2}{2})}} \sin^{s-3}(\pi t_2) &\text{for } 0 \le t_2 \le 1 \\
0 &\text{otherwise.}
\end{cases}
\end{equation*}
Let $\gamma = (1/\sigma) V_{22}^{-1/2} \tau$.
Define $f_Q$ to be the probability density function of a noncentral chi squared distribution
with $s$ degrees of freedom and noncentrality parameter $||\gamma||^2$.
Also define
\begin{equation*}
d \big(t_1, r; s, ||\gamma|| \big) =
\begin{cases}
||\gamma||^2 + 2 ||\gamma||r \cos(2 \pi t_1) + r^2 & \text{for } s = 2 \\
||\gamma||^2 + 2 ||\gamma||r \cos(\pi t_1) + r^2 & \text{for } s \ge 3.
\end{cases}
\end{equation*}
Define the unit vector $u_b = (1/||b||) b$. When $||\gamma|| > 0$, define $u_{\gamma} = (1/||\gamma||) \gamma$ and then
$\psi = u_b^T u_{\gamma}$. Also define $\psi = 1$ when $||\gamma|| = 0$. Now, when $||\gamma|| > 0$, define
\begin{align*}
k(t_1 ; \psi) &=
\psi \cos(2 \pi t_1) + \sqrt{1 - \psi^2} \, \sin(2 \pi t_1) \\
k(t_1, t_2; s, \psi) &=
\begin{cases}
\psi \cos(\pi t_1) + \sqrt{1 - \psi^2} \, \sin(\pi t_1) \cos(2 \pi t_2) & \text{for } s = 3 \\
\psi \cos(\pi t_1) + \sqrt{1 - \psi^2} \, \sin(\pi t_1) \cos(\pi t_2) & \text{for } s \ge 4.
\end{cases}
\end{align*}
The following is the main result of the paper.

\begin{theorem}
The coverage probability of the naive $1-\alpha$ confidence interval for $\theta$ is
$P(\theta \in I, F > \ell) + P(\theta \in J, F \le \ell)$. A computationally-convenient
expression for the second term in this sum is
\begin{equation}
\label{A_J}
P(\theta \in J, F \le \ell) = \int_0^{\infty} \int_0^{s \ell w^2}
j \big(||b|| \, ||\gamma|| \, \psi, q, w; ||b|| \big) \, f_Q(q) \, f_W(w) \, dq \, dw
\end{equation}
and computationally-convenient expressions for $P(\theta \in I, F > \ell)$
are as follows.
Let $u = \sqrt{d(t_1,r;s,||\gamma||)/s \ell}$. For $s=2$,
\begin{equation}
\label{A_I_gamma_ne_0_s_2}
P(\theta \in I, F > \ell) = \int_0^1 \int_0^{\infty} \int_0^{u}
i \big(-||b|| \, r \, k(t_1; \psi), w; ||b|| \big)\,
f_W(w)\, f_R(r) \, dw \, dr \,  dt_1
\end{equation}
For $s \ge 3$ and $||\gamma||>0$,
$P(\theta \in I, F > \ell)$ is equal to
\begin{equation}
\label{quadruple}
\int_0^1 \int_0^1 \int_0^{\infty} \int_0^{u}
i \big(-||b|| \, r \, k(t_1, t_2; s, \psi), w; ||b|| \big)\,
f_W(w)\, f_R(r) \, f_{T_1}(t_1)  \, f_{T_2}(t_2) \, dw \, dr \, dt_1 \, dt_2
\end{equation}
For $s \ge 3$, $||\gamma||>0$ and $\psi \in \{-1, 1\}$,
\begin{equation*}
P(\theta \in I, F > \ell) = \int_0^1 \int_0^{\infty} \int_0^{u}
i \big(-||b|| \, r \, cos(\pi t_1), w; ||b|| \big)\,
f_W(w)\, f_R(r) \, f_{T_1}(t_1)  \, dw \, dr \, dt_1
\end{equation*}
For $s \ge 3$ and $||\gamma||=0$,
\begin{equation*}
P(\theta \in I, F > \ell) = \int_0^1 \int_0^{\infty} \int_0^{r/\sqrt{s \ell}}
i \big(-||b|| \, r \, \cos(2 \pi t_1), w; ||b|| \big)\,
f_W(w)\, f_R(r) \, dw \, dr \,  dt_1
\end{equation*}
Note that for given $||b||$ (which is determined by $a$, $C$ and $X$) and $m$, $s$, $\ell$ and $\alpha$,
the coverage probability of the naive $1-\alpha$ confidence interval is a function of
$\big ( ||\gamma||,\psi \big)$.
\end{theorem}

\medskip

\noindent The proof of this theorem is presented in Appendix A.

The formulas given in Theorem 1 have three attractive features. The first of these is that, irrespective of
how large $s$ is, these formulas involve, at most, a 4-dimensional integral.
The second of these features is that
the numerical evaluation of these integrals, reviewed in Appendix B, is very straightforward.
Thirdly, for given $m$, $s$, $\ell$, $\alpha$ and $||b||$, the coverage probability of the naive confidence
interval
is a function of the two-dimensional parameter vector $\big ( ||\gamma||,\psi \big)$. These features make
it is easy to compute the minimum coverage probability of the naive $1-\alpha$ confidence interval for
given $m$, $s$, $\ell$, $\alpha$ and $||b||$.
Finally, Theorem 1 has the following corollary
\begin{corollary}
For given $m$, $s$, $\ell$ and $\alpha$,
the infimum over $\big ( ||\gamma||,\psi \big)$ of the coverage probability of the naive $1-\alpha$ confidence interval
for $\theta$ is a function of $||b||$.
\end{corollary}

For the numerical example described in the next section, $m=4$, $s=4$, $\ell = 6.5914$
(corresponding to a 0.05 level of significance of the preliminary F test)
and $\alpha = 0.05$.
For these values of $m$, $s$ and $\ell$, the minimum coverage
probability of the naive 0.95 confidence interval, as a function of
$||b|| \in (0, 1]$, is as shown in Figure 1. All of the computations presented in the present
paper were performed with programs written in MATLAB, using the optimization and
statistics toolboxes. We note the decrease in the minimum coverage probability
of this naive confidence interval with increasing $||b||$.
We can use Corollary 1 to delineate
general categories of $a$, $C$ and $X$ (via their relationship to $||b||$)
for which this naive 0.95 confidence interval
has poor coverage properties. Specifically, this naive confidence interval
will have poor coverage properties for those values of $a$, $C$ and $X$
such that
\begin{equation*}
||b|| = \sqrt{\frac{a^T (X^TX)^{-1} C \big ( C^T (X^TX)^{-1} C \big)^{-1} C^T (X^TX)^{-1} a}
{a^T  (X^TX)^{-1} a} }
\end{equation*}
is greater than 0.7, say.

\bigskip
\bigskip


\FloatBarrier

\begin{figure}[h]
\label{Figure1}
    \centering
    \includegraphics[scale=1.1]{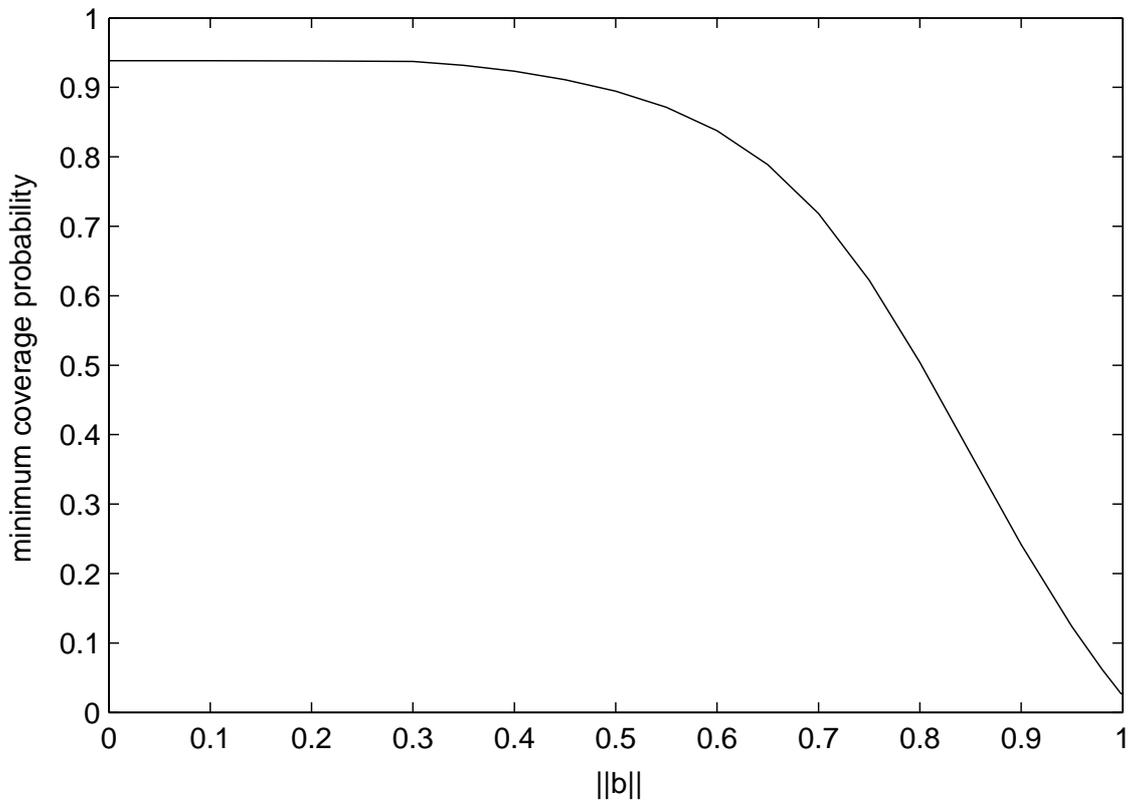}
    \caption{Plot of the minimum coverage
probability of the naive 0.95 confidence interval, as a function of
$||b|| \in (0, 1]$, for $m=4$, $s=4$, $\ell = 6.5914$
and $\alpha = 0.05$.}
\end{figure}
\FloatBarrier


\newpage

\noindent {\large \textbf{ 4. Application to a real-life data set}}

\medskip

\noindent In this section we consider the real-life {\sl analysis of covariance}
data set due to Chin et al (1994)
and analysed by Yandell (1997, Chapter 17), who makes this data available at the
website http://www.stat.wisc.edu/$\sim$yandell/pda/. This data is listed in the
Table 1. It consists of the observed response (weight gain) for a given
treatment and value of the covariate (feed intake). There are 4 possible treatments,
numbered 1 to 4.

\FloatBarrier

\begin{table}[h]
\begin{center}
\begin{tabular}{|c|c|c|}
\hline
treatment &  weight gain &  feed intake   \\
\hline
1      &     1416.1   &       2451.75 \\
1      &     1447.0   &       2546.00 \\
1      &     1509.6   &       2657.00 \\
\hline
2      &     1497.8   &       2452.10 \\
2      &     1469.9   &       2404.90 \\
2      &     1469.4   &       2479.90 \\
\hline
3      &     1510.1   &       2788.70 \\
3      &     1423.0   &       2655.50 \\
3      &     1295.9    &      2366.40 \\
\hline
4      &     1354.8   &       2578.80 \\
4      &     1326.8   &       2384.40 \\
4      &     1335.1   &       2477.50 \\
\hline
\end{tabular}
\end{center}
\caption{The observed response (weight gain) for a given
treatment and value of the covariate (feed intake).
Source: http://www.stat.wisc.edu/$\sim$yandell/pda/.}
\end{table}

\FloatBarrier

We use the following linear regression model for this data:
\begin{equation*}
Y_{ij} = \mu_i + \tilde{\beta}_i (x_{ij} - \bar{x}_{\cdot \cdot}) + \varepsilon_{ij}
\end{equation*}
where $Y_{ij}$ is the response of the $j$ th experimental unit ($j = 1, \ldots, 3$)
that is receiving the $i$ th treatment ($i=1,\ldots,4$), when the covariate takes the
value $x_{ij}$. The $\varepsilon_{ij}$ are independent and identically $N(0, \sigma^2)$
distributed and $\sigma^2$ is an unknown positive parameter.
The $\mu_i$ and  $\tilde{\beta}_i$ ($i=1,\ldots,4$) are unknown
parameters.  Also, $\bar{x}_{\cdot \cdot}$
denotes the mean of the $x_{ij}$ ($i=1,\ldots,4; \, j=1, \ldots,3$).

We express this model in the form
$Y = X \beta + \varepsilon$, where
\begin{align*}
Y &= \big(Y_{11}, Y_{12}, Y_{13},Y_{21}, Y_{22}, Y_{23},Y_{31}, Y_{32}, Y_{33},Y_{41}, Y_{42}, Y_{43} \big) \\
\varepsilon &= \big(\varepsilon_{11}, \varepsilon_{12}, \varepsilon_{13},\varepsilon_{21}, \varepsilon_{22},
\varepsilon_{23},\varepsilon_{31}, \varepsilon_{32}, \varepsilon_{33},\varepsilon_{41}, \varepsilon_{42}, \varepsilon_{43} \big),
\end{align*}
$\beta = \big(\mu_1, \mu_2, \mu_3, \mu_4, \tilde{\beta}_1, \tilde{\beta}_2, \tilde{\beta}_3, \tilde{\beta}_4 \big)$
and $X$ is the obvious $12 \times 8$ matrix.

As considered by Yandell (1997, p.271), we carry out a preliminary test of the null hypothesis
$H_0: \tilde{\beta}_1 = \tilde{\beta}_2 = \tilde{\beta}_3 = \tilde{\beta}_4$ against the alternative
hypothesis  $H_A$ that the $\tilde{\beta}_i$ are not all the same, using an F test. Suppose that we use
a 0.05 level of significance for this preliminary test. We express $H_0$ as $\tau = 0$
and $H_A$ as $\tau \ne 0$, where
$\tau = C^T \beta$ and
\begin{equation*}
C^T = \left[\begin{matrix} 0 \quad 0 \quad 0 \quad 0 \quad -1 \quad 1 \quad 0 \quad 0 \\
0 \quad 0 \quad 0 \quad 0 \quad -1 \quad 0 \quad 1 \quad 0  \\
0 \quad 0 \quad 0 \quad 0 \quad -1 \quad 0 \quad 0 \quad 1  \\
\end{matrix} \right ]
\end{equation*}

Define the parameter of interest $\theta$ as follows. Let $Y_1^*$ and $Y_2^*$ denote
the responses of two experimental units, receiving treatments 1 and 2 respectively, for
the same value $x^*$ of the covariate. In other words,
\begin{align*}
Y_1^* &= \mu_1 + \tilde{\beta}_1 ( x^* - \bar{x}_{\cdot \cdot}) + \varepsilon_1^* \\
Y_2^* &= \mu_2 + \tilde{\beta}_2 ( x^* - \bar{x}_{\cdot \cdot}) + \varepsilon_2^*
\end{align*}
where $\varepsilon_1^*$ and $\varepsilon_2^*$ are independent and identically $N(0, \sigma^2)$
distributed random variables. Let
$\theta = E(Y_1^*) - E(Y_2^*) = \mu_1 - \mu_2 + (\tilde{\beta}_1 - \tilde{\beta}_2)(x^* - \bar{x}_{\cdot \cdot})$.
Thus, $\theta = a^T \beta$ where
$a = \big(1, -1, 0, 0, (x^* - \bar{x}_{\cdot \cdot}), -(x^* - \bar{x}_{\cdot \cdot}), 0, 0 \big)$.
We suppose that $x^* - \bar{x}_{\cdot \cdot} = 125.39$, which is the maximum value of
$|x_{ij} - \bar{x}_{\cdot \cdot}|$ for the data.

For this situation, $||b|| = 0.96869$ and so the minimum coverage
probability of the naive 0.95 confidence interval is 0.0846.
This shows that this confidence interval is completely inadequate,
in this situation.

\bigskip

\noindent {\large{\bf 5. Discussion}}

\medskip

\noindent The poor coverage properties of naive confidence intervals found in this paper
are presaged by the poor coverage properties of naive confidence intervals
found in the context of a preliminary best subset variable selection by minimizing an
AIC-type criterion, see e.g. Kabaila (2005), Kabaila \& Leeb (2006) and Kabaila \& Giri
(2009). Apart from the form of preliminary model selection used, minimum AIC versus an F test,
these papers differ from the present paper in that the present paper provides a method
for computing the minimum coverage probability, whereas Kabaila (2005), Kabaila \& Leeb (2006) and Kabaila \& Giri
(2009) provide only upper bounds on the minimum coverage probability of the naive
confidence interval.


\newpage

\noindent {\large{\bf Appendix A: Proof of Theorem 1}}

\medskip

\noindent In this appendix, we prove Theorem 1.
Define $G = (\hat \Theta - \theta)/(\sigma \sqrt{v_{11}})$ and $H = (1/\sigma) V_{22}^{-1/2} \hat \tau$.
Let $f_H$ denote the probability density function of $H$.
Note that
\begin{equation}
\label{model_G_H}
\left[\begin{matrix} G\\ H \end{matrix}
\right] \sim \text{N} \left ( \left[\begin{matrix} 0 \\ \gamma \end{matrix}
\right], \left[\begin{matrix} 1 \quad b^T \\ b \quad I_s \end{matrix}
\right] \right ).
\end{equation}
where $I_s$ denotes the $s \times s$ identity matrix. Thus the distribution of $G$, conditional
on $H = h$, is $N \big( b^T(h-\gamma), 1 - ||b||^2 \big)$. Note that $(G,H^T)^T$ and $W$ are
independent random vectors.
We use the notation
\begin{equation*}
\chi({\cal A}) =
\begin{cases}
1 &\text{if } {\cal A} \ \ \text{is true} \\
0 &\text{if } {\cal A} \ \ \text{is false}
\end{cases}
\end{equation*}
where ${\cal A}$ is an arbitrary statement. This is similar to the Iverson bracket notation
(Knuth, 1992).

By the law of total probability, the coverage probability of the naive $1-\alpha$ confidence interval is
$$P(\theta \in I, F > \ell) + P(\theta \in J, F \le \ell).$$
We divide the remainder of the proof into 2 parts.

\medskip

\noindent {\bf Part 1: expression for $\boldsymbol {P(\theta \in I, F > \ell)}$}

\smallskip

\noindent Suppose that $||\gamma||>0$. We prove the validity of the expressions
\eqref{A_I_gamma_ne_0_s_2} and \eqref{quadruple} for $P(\theta \in I, F > \ell)$.
The proofs of the validity of the other expressions for $P(\theta \in I, F > \ell)$
(given in Theorem 1) are similar and are omitted, for the sake of brevity.

Now
$\{ \theta \in I \} = \{ -t(m) W \le G \le t(m) W \}$ and $F = H^T H/(s W^2)$.
Thus
\begin{align}
&P(\theta \in I, F > \ell) \notag \\
&= P \left ( -t(m) W \le G \le t(m) W, \frac{H^T H}{s W^2} > \ell \right)  \notag \\
&= \int \cdots \int \int_0^{\infty}
P \left(-t(m) W \le G \le t(m) W, \frac{H^T H}{s W^2} > \ell \, \Big| \, W=w, \, H=h \right) \, f_W(w) \, dw \, f_H(h) \, dh \notag \\
\label{first_interm}
&= \int \cdots \int \int_0^{\infty}
P \left(-t(m) w \le G \le t(m) w \, \Big| \, H=h \right) \, \chi \big(h^T h > s \, \ell \, w^2 \big) \, f_W(w) \, dw \, f_H(h) \, dh
\end{align}
Note that
\begin{align*}
&P (-t(m) w \le G \le t(m) w \, | \, H=h) \\
&= P \big(-t(m) w - b^T (h - \gamma) \le Z \le t(m) w - b^T (h - \gamma) \big) \quad \text{where } Z \sim N(0, 1 - ||b||^2) \\
&= i\big(-b^T (h - \gamma), w; ||b|| \big).
\end{align*}
Thus
\begin{align}
\label{second_interm}
\eqref{first_interm} &= \int \cdots \int \int_0^{\infty}
i\big(-b^T (h - \gamma), w; ||b|| \big) \, \chi \big(h^T h > s \, \ell \, w^2 \big) \, f_W(w) \, dw \, f_H(h) \, dh
\notag \\
&= E \Big( i \big(-b^T (H-\gamma), W; ||b|| \big) \, \chi \big(H^T H > s\, \ell \, W^2 \big) \Big) \notag \\
&= E \Big( i \big(-b^T (H-\gamma), W; ||b|| \big) \, \chi \big(W \le \sqrt{H^T H / s \ell}  \big) \Big).
\end{align}

We now find a simple formula for this expected value.
Since $H \sim N(\gamma, I_s)$, $H = \gamma + R U$ where $R$ is a nonnegative random variable and $U$ is a random
$s$-vector with the following distributions. The random vectors $R^2$ and $U$ are independent, with
$R^2 \sim \chi^2_s$ and $U$ is distributed uniformly on the surface of the unit sphere
in $\mathbb{R}^s$. Thus $b^T (H - \gamma) = R \, ||b|| \, L_b$, where $L_b = u_b^T U$.
Also, $H^T H = ||\gamma||^2 + 2 ||\gamma|| R L_{\gamma} + R^2$, where $L_{\gamma} = u_{\gamma}^T U$.
 Hence
 \begin{equation*}
 \eqref{second_interm} =
 E \Big( i \big(-||b|| R L_b, W; ||b|| \big) \,
 \chi \big(W \le \sqrt{(||\gamma||^2 + 2 ||\gamma|| R L_{\gamma} + R^2) / s\ell} \big) \Big).
  \end{equation*}
 Note that $(L_{\gamma}, L_b)$, $R$ and $W$ are independent random vectors.
 Define the random vector $(T_1, T_2)$
 to be such that $T_1$, $T_2$, $R$ and $W$ are independent and $T_1$ and $T_2$ have the probability density functions
 $f_{T_1}$  and $f_{T_2}$ respectively, defined in Section 3.
 Define the unit $s$-vectors $e_{\gamma}$ and $e_b$ as follows.
 The vector $e_{\gamma}$ has 1 as its first component and zeros for the remaining components.
 The vector $e_b$ has first component $\psi$, second component $\sqrt{1 - \psi^2}$
 and zeros for the remaining components.
Because $U$ is distributed uniformly on the surface of the unit sphere
in $\mathbb{R}^s$, $(L_{\gamma}, L_b)$ has the same distribution as $(e_{\gamma}^T U, e_b^T U)$.
It follows from Fang and Wang (1994, p.49, pp.306--306 and p.308) that $(U_1, U_2)$ has the same distribution as
(a) $\big ( \cos (2 \pi T_1), \, \sin (2 \pi T_1) \big)$ for $s=2$, (b)
$\big ( \cos (\pi T_1), \, \sin (\pi T_1) \cos (2 \pi T_2) \big)$ for $s=3$ and
$\big ( \cos (\pi T_1), \, \sin (\pi T_1) \cos (\pi T_2) \big)$ for $s>3$.
Thus, $(L_{\gamma}, L_b)$ has the same distribution as
\begin{equation*}
\Big(\cos(2 \pi T_1), \;
\psi \cos(2 \pi T_1) + \sqrt{1 - \psi^2} \sin (2 \pi T_1) \Big)
\end{equation*}
for $s=2$,
\begin{equation*}
\Big(\cos(\pi T_1), \;
\psi \cos(\pi T_1) + \sqrt{1 - \psi^2} \sin (\pi T_1) \cos (2 \pi T_2) \Big).
\end{equation*}
for $s=3$ and
\begin{equation*}
\Big(\cos(\pi T_1), \;
\psi \cos(\pi T_1) + \sqrt{1 - \psi^2} \sin (\pi T_1) \cos (\pi T_2) \Big).
\end{equation*}
for $s>3$.
Thus $P(\theta \in I, F > \ell)$ is
 \begin{equation*}
 E \Big( i \big(-||b|| \, R \, k(T_1;\psi), W; ||b|| \big) \,
 \chi \big(W \le \sqrt{d(T_1, R; 2, ||\gamma||) / s \ell} \big) \Big)
 \end{equation*}
 for $s=2$ and
 \begin{equation}
 E \Big( i \big(-||b|| \, R \, k(T_1,T_2;s,\psi), W; ||b|| \big) \,
 \chi \big(W \le \sqrt{(d(T_1, R; s, ||\gamma||) / s \ell} \big) \Big)
 \end{equation}
 for $s \ge 3$.
 This leads to the expressions \eqref{A_I_gamma_ne_0_s_2} and \eqref{quadruple}
for $P(\theta \in I, F > \ell)$ given in the theorem.

 \medskip

\noindent {\bf Part 2: expression for $\boldsymbol {P(\theta \in J, F \le \ell)}$}

\smallskip

\noindent The derivation of the expression for $P(\theta \in J, F \le \ell)$
is based on \eqref{J} and the fact that $\Theta^*$, $\hat \tau$ and $R(\hat \beta)$
are independent random vectors. Define $Q = (1/\sigma^2) \hat \tau^T V_{22}^{-1} \hat \tau$
and note that $Q$ has a noncentral chi squared distribution with $s$ degrees of freedom and
noncentrality parameter $||\gamma||^2$. Note that
\begin{equation*}
\{ \theta \in J \} =
\left \{ Z \in \left [ b^T \gamma \pm t(m+s) \sqrt{\frac{m W^2 + Q}{m+s}} \sqrt{1 - ||b||^2} \right ] \right \}
\end{equation*}
where $Z$, $Q$ and $W$ are independent random variables and $Z \sim N(0, 1 - ||b||^2)$. Also
\begin{equation*}
F = \frac{\hat \tau^T V_{22}^{-1} \hat \tau / s}{\hat \Sigma^2} = \frac{Q/s}{W^2}.
\end{equation*}
Thus
\begin{align*}
P(\theta \in J, F \le \ell)
&= P \left ( Z \in \left [ b^T \gamma \pm t(m+s) \sqrt{\frac{m W^2 + Q}{m+s}} \sqrt{1 - ||b||^2} \right ], \,
Q \le s \ell W^2 \right ) \\
&= \int_0^{\infty} \int_0^{s \ell w^2}
j \big(||b|| \, ||\gamma|| \, \psi, q, w; ||b|| \big) \, f_Q(q) \, f_W(w) \, dq \, dw
\end{align*}
by a method similar to that used in Part 1.

\bigskip

\noindent {\large{\bf Appendix B: Numerical evaluation of the integrals in \newline Theorem 1}}

\medskip

\noindent We evaluate the integrals
\eqref{A_J} and \eqref{A_I_gamma_ne_0_s_2}
in the statement of Theorem 1 as follows.
We approximate \eqref{A_J} by
\begin{equation}
\label{approx_A_J}
\int_0^{c_1} \int_0^{s \ell w^2}
j \big(||b|| \, ||\gamma|| \, \psi, q, w; ||b|| \big) \, f_Q(q) \, f_W(w) \, dq \, dw
\end{equation}
for an appropriately chosen value of $c_1$.
We bound the error of this approximation as follows.
Since $j \big(||b|| \, ||\gamma|| \, \psi, q, w; ||b|| \big)$ is a probability,
\begin{align*}
0 \le \eqref{A_J} - \eqref{approx_A_J}
&\le \int_{c_1}^{\infty} \int_0^{s \ell w^2}
 \, f_Q(q) \, f_W(w) \, dq \, dw \\
&\le \int_{c_1}^{\infty}\, f_W(w)\, dw \\
&= P(M > m \, c_1^2) \quad \text{where }\ M \sim \chi_m^2.
\end{align*}
We choose $c_1$ sufficiently large that the right hand side is less than, say, $10^{-5}$.

To evaluate \eqref{approx_A_J}, we transform the region of integration to a rectangle as follows.
Change the variable of integration $q$ in \eqref{approx_A_J} to $q^* = q/(s \ell w^2)$,
so that \eqref{approx_A_J} is equal to
\begin{equation*}
\int_0^{c_1} \int_0^1
j \big(||b|| \, ||\gamma|| \, \psi, s \ell w^2 q^*, w; ||b|| \big) \, f_Q(s \ell w^2 q^*) \,
s \ell w^2 \, f_W(w) \, dq^* \, dw.
\end{equation*}
The integrand is a smooth function of $(q^*, w) \in [0,1] \times [0, c_1]$
and so this double integral is easily evaluated by numerical integration.

We approximate \eqref{A_I_gamma_ne_0_s_2}
by
\begin{equation}
\label{approx_A_I}
\int_0^1 \int_0^{c_2} \int_0^{u}
i \big(-||b|| \, r \, k(t_1; \psi), w; ||b|| \big)\,
f_W(w)\, f_R(r) \, dw \, dr \,  dt_1
\end{equation}
for an appropriately chosen value of $c_2$.
We bound the error of this approximation as follows.
Since $i \big(-||b|| \, r \, k(t_1; \psi), w; ||b|| \big)$ is a probability,
\begin{align*}
0 \le \eqref{A_I_gamma_ne_0_s_2} - \eqref{approx_A_I}
&\le \int_0^1 \int_{c_2}^{\infty} \int_0^u  f_W(w)\, dw \, f_R(r)\, dr \, f_T(t) \, dt \\
&\le \int_{c_2}^{\infty}\, f_R(r)\, dr \\
&= P(R^2 > c_2^2)\quad \text{where }\ R^2 \sim \chi_s^2.
\end{align*}
We choose $c_2$ sufficiently large that the right hand side is less than, say, $10^{-5}$.

To evaluate \eqref{approx_A_I}, we transform the region of integration to a rectangle as follows.
Change the variable of integration $w$ in \eqref{approx_A_I} to $w^* = w/\sqrt{d(t_1,r;s,||\gamma||)/s \ell}$,
so that \eqref{approx_A_I} is equal to
\begin{align*}
&\int_0^1 \int_0^{c_2} \int_0^1
i \big(-||b|| \, r \, k(t_1; \psi), \sqrt{d(t_1,r;s,||\gamma||)/s \ell} \, w^*; ||b|| \big)\,
f_W \big(\sqrt{d(t_1,r;s,||\gamma||)/s \ell} \, w^* \big)\\
&\phantom{1234567890123456789012345678901234567890}
\sqrt{d(t_1,r;s,||\gamma||)/s \ell} \, f_R(r) \, dw^* \, dr \,  dt_1.
\end{align*}
The integrand is a smooth function of $(w^*, r, t_1) \in [0,1] \times [0, c_2] \times [0,1]$
and so this triple integral is easily evaluated by numerical integration.

The evaluation of the other integrals in the statement of Theorem 1 is similar to the evaluation
of the integrals \eqref{A_J} and \eqref{A_I_gamma_ne_0_s_2}.
The evaluation of \eqref{quadruple} using MATLAB requires special
comment. In MATLAB, the highest dimensional integral that one can evaluate using a built-in MATLAB
function is a triple integral. We evaluate the quadruple integral \eqref{quadruple} using MATLAB
as follows. As before, let $u = \sqrt{d(t_1,r;s,||\gamma||)/s \ell}$. Define
\begin{equation*}
g(w^*, t_1, r) = \int_0^1
i \big(-||b|| \, r \, k(t_1, t_2; s, \psi), \, u \, w^*; ||b|| \big)
\, f_{T_2}(t_2) \, dt_2
\end{equation*}
The integrand on the right-hand-side is a very smooth function of $t_2$. We evaluate $g(w^*, t_1, r)$,
to a good approximation, using a compound Simpson's rule with a specified number of subdivisions
of the interval $[0,1]$. The quadruple integral \eqref{quadruple} is approximated by
\begin{equation*}
\int_0^{c_2} \int_0^1 \int_0^1
g(w^*, t_1, r) \,
f_W(u \, w^*)\, u \,
f_{T_1}(t_1) \, f_R(r) \,  dw^*  \, dt_1 \, dr
\end{equation*}
which is evaluated using the MATLAB built-in function triplequad.

\bigskip

\nin{\large{\bf References}}

\medskip

\rf Chin, S.F., Storkson, J.M., Albright, K.J., Cook, M.E. \& Pariza, M.W.:  Conjugate linoleic
acid is a growth factor for rats as shown by enhanced weight gain and improved feed efficiency.
Journal of Nutrition 124, 2344 -- 2349 (1994)

\smallskip

\rf Fang, K.-T. \& Wang, Y.:  Number-theoretic Methods in Statistics.
Chapman \& Hall, London (1994)

\smallskip

\rf Farchione, D.: Interval estimators that utilize uncertain prior information. Unpublished Ph.D. thesis,
Department of Mathematics and Statistics, La Trobe University (2009)

\smallskip

\rf Freund, R.J., Wilson, W.J. \& Sa, P.:  Regression Analysis: Statistical Modeling
of a Response Variable, 2nd ed.. Elsevier, Academic Press, Burlington, Mass. (2006)

\smallskip

\rf Graybill, F.A.:  Theory and Application of the Linear Model. Duxbury, Pacific Grove, CA (1976)

\smallskip

\rf Kabaila, P.:  On the coverage probability of confidence intervals
in regression after variable selection.  Australian \& New Zealand Journal of Statistics
47, 549--562 (2005).

\smallskip

\rf Kabaila, P., Leeb, H.:  On the large-sample minimal coverage
probability of confidence intervals after
model selection. Journal of the American Statistical Association 101, 619--629 (2006)

\smallskip

\rf Kabaila, P., Giri, K.:  Upper bounds on the minimum coverage probability of confidence
intervals in regression after model selection. Australian \& New Zealand Journal of Statistics
 51, 271 -- 288 (2009)

\smallskip

\rf Knuth, D.E.:  Two notes on notation. American Mathematical Monthly  99, 403--422 (1992)

\smallskip

\rf Kuehl, R.O.:  Design of Experiments: Statistical Principles of Research Design and Analysis,
2nd ed.. Brooks/Cole,Pacific Grove, CA (2002)

\smallskip

\rf Milliken, G.A., Johnson, D.E.:  Analysis of Messy Data, Volume III: Analysis of
Covariance. Chapman \& Hall/CRC, Boca Raton, Fl. (2002)

\smallskip

\rf Yandell, B.S.:  Practical Data Analysis for Designed Experiments.
Chapman \& Hall, London, New York (1997).

\end{document}